\documentclass[]{article}
\usepackage[T1]{fontenc}
\usepackage[utf8]{inputenc}
\usepackage{authblk}
\usepackage{amsmath,amssymb,bm}
\usepackage{epsfig,graphicx,subfigure,color}
\usepackage{multirow,booktabs,fancyhdr,setspace}
\usepackage{longtable}  
\usepackage{hyperref}   
\usepackage{array}
\usepackage{natbib}
\usepackage{enumitem}
\usepackage{dsfont}

\graphicspath{{figs/}}

\newtheorem{theorem}{Theorem}
\newtheorem{lemma}{Lemma}

\def\MSE{\textsc{mse}}

\def\EFF{\textsc{eff}}

\begin{document}
\title{%
Optimal group testing designs for estimating prevalence with uncertain testing errors \footnote{{\it Address for correspondence}:
Mong-Na Lo Huang,
Department of Applied Mathematics, National Sun Yat-sen University,
Kaohsiung 804, Taiwan. 
E-mail: lomn@math.nsysu.edu.tw}
}
\author[1,2,3]{Shih-Hao Huang}
\author[1]{Mong-Na Lo Huang}
\author[3]{Kerby Shedden}
\author[4]{Weng Kee Wong}
\affil[1]{National Sun Yat-sen University, Kaohsiung, Taiwan.}
\affil[2]{Academia Sinica, Taipei, Taiwan.}
\affil[3]{University of Michigan, Ann Arbor, U.S.A.}
\affil[4]{University of California at Los Angeles, U.S.A.}

\date{}

\maketitle

\begin{abstract}
We construct optimal designs for group testing experiments 
where the goal is to estimate the prevalence of a trait 
by using a test with uncertain sensitivity and specificity. 
Using optimal design theory for approximate designs, 
we show that the most efficient design 
for simultaneously estimating the prevalence, sensitivity and specificity 
requires three different group sizes with equal frequencies. 
However, if estimating prevalence as accurately as possible is the only focus, 
the optimal strategy is to have three group sizes with unequal frequencies. 
On the basis of a chlamydia study in the U.S.A., 
we compare performances of competing designs and
provide insights into how the unknown sensitivity and specificity of the test 
affect the performance of the prevalence estimator.
We demonstrate that the locally $D$- and $D_s$-optimal designs
proposed have high efficiencies even when 
the prespecified values of the parameters are moderately misspecified.

\vspace{2mm}
\end{abstract}
%

\section{Introduction}

In a group testing study, 
the goal is often to estimate the prevalence of a rare disease or a particular trait 
\citep{HOS94,HOR00}. 
Group testing is frequently used in studies 
where testing individuals for a trait is expensive and
individual samples are relatively plentiful. 
In a group testing study, samples from individuals are pooled and tested as a single unit. 
Fewer tests or trials are therefore required so that the
cost of the study is reduced. 
(Throughout, we use the terms trials and tests interchangeably.)
The key assumption is that the test result from a group of individuals is positive if and only if
at least one individual in the group has the trait.

Applications of group testing abound in the literature, 
especially for case diagnosis \citep{McMahan12} and 
prevalence estimation of diseases \citep{Shipitsyna07}. 
Several researchers have shown that group testing techniques 
can be used to estimate the prevalence of a rare disease efficiently \citep{HO06}. 
\cite{Dorf43} gave an early perspective and an overview. 
From the design point of view, some key research questions are to determine 
how many distinct group sizes should be used, 
and how many trials at each such group size should be run. 
In practice, group testing studies usually have only one group size, 
but this may be because the potential advantages of using unequal group sizes 
have not been well explored.

Much of the work on group testing makes a simplifying assumption that 
there is no testing error; see for example, \cite{HOR00}. 
In practice, the test in most group testing experiments is likely to have non-zero error rates, 
with sensitivity and specificity being less than $100\%$. 
We recall that sensitivity is the probability that 
the test will be positive given that the sample tested is truly positive, 
and specificity is the probability that 
the test will be negative given that the sample tested is truly negative.
A few researchers have considered application of group testing 
when there may be testing errors, 
but they assumed that the sensitivity and specificity 
are known and do not depend on the group size.
When the total number of individuals is fixed and all group sizes are equal, 
\cite{Tu95} obtained the maximum likelihood estimator (MLE) of the prevalence 
and determined graphically the common group size that 
minimizes the variance of the prevalence estimator. 
For a given group size, \cite{Liu12} determined conditions under which group testing performs
better than individual testing, when either the total number of trials or the total number of
individuals is fixed.
The sensitivity and specificity can be directly estimated 
if samples are tested by using both the test of interest 
and a ``gold standard'' test that has no testing errors. 
Such a gold standard test may not exist, however, and, even if one does exist, 
it may be too expensive or complex for routine use. 
Lacking a gold standard test, investigators may 
prespecify sensitivity and specificity values.
However, nominal values may be optimistic, 
and values that are estimated from different populations could be biased. 
Therefore it is more realistic in applications to assume that 
the sensitivity and specificity of the test are both unknown. 
Under this assumption, \cite{Zhang14} used numerical methods 
to obtain group testing designs with or without a gold standard test. 
Their proposed designs are not fully supported by theory as being optimal.
When the cost of a test is greater than the cost of collecting an individual sample, 
it may be reasonable to design the study by taking the total number of trials, 
rather than the total number of individuals, to be fixed. 
Our goal is to apply optimal design theory to construct optimal group testing designs 
when we have a fixed number of trials, and 
when both the sensitivity and the specificity of the test are unknown. 
We consider two situations:
(a) we want to estimate the prevalence, sensitivity and specificity, all with equal interest;
(b) only the prevalence is of interest, 
but both testing error rates are treated as nuisance parameters in the design problem.
In the framework of optimal experimental design \citep{atk07,wong09}, 
the first setting requires estimating all unknown parameters, 
and the second setting requires the estimation of 
only one of the three unknown parameters in the model.
$D$-optimality as a design criterion is appropriate for the first situation 
and $D_s$-optimality is appropriate for the second situation. 
The resulting optimal designs are called $D$-optimal 
and $D_s$-optimal designs respectively. 
These optimal designs minimize the generalized variance of 
all or a subset of the model parameters of interest among all designs, 
so they provide the most accurate inference for all three or just the prevalence parameter. 
Other design criteria, such as $A$-optimality or $c$-optimality, can also be used.
Our set-up is different from earlier studies for 
finding efficient or optimal designs for group testing experiments, 
such as \cite{Tu95} and \cite{Zhang14}. 
The main differences are that
(a) our optimal designs are sought among all possible designs 
without any restriction on the number of different group sizes, and
(b) our optimal designs are justified by theory.
Tu?s results are more restrictive because their designs are sought 
among all designs with a common group size and 
cannot incorporate uncertainty in testing error rates. 
Zhang?s designs are obtained by numerical optimization 
with dimension as large as the upper bound of the group sizes; 
therefore their global optimality cannot be guaranteed, 
and we may not be able to interpret 
how prespecified parameter values affect their optimal designs. 
To our best knowledge,
our work is the first that uses optimal design theory 
to provide a theoretical justification for designing group testing experiments 
when testing error rates are uncertain.
The organization of our paper is as follows. 
Section 2 describes the statistical background.
Section 3 constructs optimal group testing designs 
for the two situations that were described above. 
In Section 4, we evaluate the finite sample properties 
and robustness properties of our proposed optimal designs
when the parameter vector is misspecified. 
The set-up of the simulation is based on a chlamydia study conducted in the U.S.A. 
We end the paper with a discussion in Section 5 
with remarks on alternative methods and 
ideas for further research directions in this area.

All computations in this paper were done by using our code 
called {\it gtest.nb} developed by using Mathematica 10.0; 
see \cite{Wolfram15} for more details. 
The package can be freely downloaded from 
\url{http://www.math.nsysu.edu.tw/~lomn/huangsh/}.

\section{Preliminaries}

Throughout, we denote the prevalence, sensitivity and specificity of the test by $p_0$, $p_1$ and $p_2$, respectively. We assume that $p_0\in(0,1)$, $p_1,p_2\in(0.5,1]$, and the two types of testing errors occur randomly with testing error rates $1-p_1$ and $1-p_2$ respectively. Letting $\theta=(p_0,p_1,p_2)^\mathrm{T}$, a direct calculation shows that the probability of a positive response (either a true positive or a false positive) for a trial with group size $x$ is
\begin{align}
\pi(x)=\pi(x|\theta)&=p_1\left(1-\left(1-p_0\right)^x\right)+\left(1-p_2\right)\left(1-p_0\right)^x
=p_1-\left(p_1+p_2-1\right)\left(1-p_0\right)^x.
\label{eq:gt_resp}
\end{align}
In practice, there are generally limits on the group sizes so we consider designs that are subject to a known group size constraint $1\leq x_L\leq x\leq x_U<\infty$. We note that $\pi(x)$ is a convex combination of $p_1$ and $1-p_2$. For a study consisting of $n$ trials at $k$ distinct group sizes $x_1<x_2\cdots<x_k$, the log-likelihood function of $\theta$ is
\begin{align}
L(\theta)=\sum_{i=1}^k\{y_i\log(\pi(x_i|\theta))+(n_i-y_i)\log(1-\pi(x_i|\theta))\},
\label{eq:loglike}
\end{align}
apart from an unimportant additive constant.  Here each $y_i$ is the number of positive responses among the $n_i$ trials at group size $x_i$ and $\sum_{i=1}^k n_i=n$. We assume that each $y_i$ has a binomial distribution with parameters $(n_i,\pi(x_i|\theta))$, and the $\{y_i\}_{i=1}^k$ are mutually independent.
In what follows, we work with approximate designs advocated by \cite{Kiefer74}. An approximate group testing design denoted by $\xi=\{(x_i,w_i)\}_{i=1}^k$ contains $k$ distinct group sizes (or support points) $x_1<x_2<\cdots<x_k$ and corresponding weights $w_1,w_2,\ldots,w_k$, where $\sum_{i=1}^kw_i=1$ and each \mbox{$w_i>0$} is interpreted as the proportion of trials at group size $x_i$. When the total number of trials $n$ is fixed a priori by either cost or practical considerations, the group testing design $\xi$ takes $nw_i$ groups of size $x_i,i=1,\ldots,k$ subject to $nw_1+\ldots+nw_k=n$. In reality, the group sizes and the numbers of trials at these sizes must be positive integers. Therefore, optimal approximate group testing designs are rounded before implementation so that each group size and its frequency are positive integers. We call such designs rounded optimal designs. Simple rounding procedures to obtain a design for implementation have been shown to produce little loss in efficiency. More efficient rounding apportionment introduced in section 12.5 of \cite{puk06} can also be applied to obtain the optimal number of trials at the rounded group size.
Under model (\ref{eq:gt_resp}), the maximum likelihood estimator (MLE) of $\theta$ is obtained by maximizing (\ref{eq:loglike}).  If design $\xi$ were used to obtain the data, a direct calculation shows the covariance matrix of the estimated $\theta$ is asymptotically proportional to the inverse of the $3\times3$ Fisher information matrix given by
\begin{align}
M(\xi)=\sum_{i=1}^k w_i\lambda(x_i)f(x_i)f(x_i)^\mathrm{T}=M_f(\xi)\cdot \textrm{Diag}_\lambda(\xi)\cdot M_f(\xi)^\mathrm{T}.\label{eq:inforM2}
\end{align}
Here
\begin{align}
\lambda(x)&=[\pi(x)(1-\pi(x))]^{-1},\nonumber\\
f(x)&=(f_0(x),f_1(x),f_2(x))^\mathrm{T},\nonumber\\
f_0(x)&=\partial\pi(x)/\partial p_0=x(p_1+p_2-1)(1-p_0)^{x-1},\nonumber\\
f_1(x)&=\partial\pi(x)/\partial p_1=1-(1-p_0)^x,\nonumber\\
f_2(x)&=\partial\pi(x)/\partial p_2=-(1-p_0)^x,\nonumber\\
M_f(\xi)&=(f(x_1),~f(x_2),~\cdots,~f(x_k))\in\mathbb{R}^{3\times k} \label{eq:Mf},
\end{align}
and
$
\textrm{Diag}_\lambda(\xi)~\text{ is a $k\times k$ diagonal matrix with elements } \{w_i\lambda(x_i)\}_{i=1}^k.
$
The design we seek for estimating $\theta$ is a $D$-optimal design that maximizes the determinant of the information matrix (\ref{eq:inforM2}), or equivalently a design that minimizes the generalized variance of all parameter estimators, by choice of $k,x_1,\ldots,x_k,w_1,\ldots,w_k$. In other words, a $D$-optimal design maximizes
\begin{align}
\Phi_D[M(\xi)]=\log|M(\xi)| \label{eq:Dcriterion}
\end{align}
among all possible designs under which $\theta$ is estimable. When we are interested in estimating only a subset of the model parameters, we use the $D_s$-optimality criterion. Here we are interested in accurately estimating only the prevalence, and we treat sensitivity and specificity as nuisance parameters. The $D_s$-criterion here is to minimize only the asymptotic variance of the prevalence estimator. Such a $D_s$-optimal design maximizes
\begin{align}
\Phi_s[M(\xi)]=-\log\left(M(\xi)^{-}\right)_{11} \label{eq:Dscriterion}
\end{align}
among all designs under which $p_0$ is estimable, where $M^{-}$ is a generalized inverse of matrix $M$ and we use the notation $O_{11}$ to denote the $(1,1)$ entry of the matrix $O$. Note that the $D_s$-optimality above is equivalent to $c$-optimality with $c=(1,0,0)^{\rm T}$ \citep{Fed72}.
The solutions to our optimization problems require solving the equations in Lemma \ref{lemma:eqs} below. Its proof mainly requires calculus and the intermediate value theorem and has therefore been omitted. Before stating Lemma \ref{lemma:eqs}, we introduce further notation. Let 
\begin{align*}
c&=p_1/\{(p_1+p_2-1)(1-p_0)^{x_L}\}>1,\\
\delta&=(1-p_1)/p_1\in[0,1),\\
r&=(1-p_0)^{x_U-x_L}\in(0,1),\quad\text{and}\\
\Delta_0&=r\log(r)/(1-r).
\end{align*}
For $a\in(r,1)$, let
\begin{align*}
\Delta_1(a)
=\frac{(a-r)\sqrt{(c-1)(\delta c+1)}+(1-a)\sqrt{(c-r)(\delta c+r)}}{(1-r)\sqrt{(c-a)(\delta c+a)}},
\end{align*}
and let 
\begin{align*}
\Delta_2(a)
=\frac{\sqrt{(c-1)(\delta c+1)}-\sqrt{(c-r)(\delta c+r)}}{(1-r)\sqrt{(c-a)(\delta c+a)}}.
\end{align*}
These definitions imply that $r<a<1<c$ which is used in the proofs without further mention.
\begin{lemma}
\label{lemma:eqs}
Each of the following equations has a root $a\in(r,1)$:
\begin{align}
&\frac{2}{a}\left(1+\frac{1+\Delta_0\times\frac{1}{a}}
{\log(a)-\Delta_0\times\left(\frac{1}{a}-1\right)}\right)
=\frac{1}{\delta c+a}-\frac{1}{c-a},\label{eq:gtmDx2}
\end{align}
and
\begin{align}
&\frac{2\left(1+\Delta_1(a)\right)}{a}\left(1+\frac{1+\Delta_0\times\frac{1}{a}}
{\log(a)-\Delta_0\times\left(\frac{1}{a}-1\right)}\right)=\frac{1}{\delta c+a}-\frac{1}{c-a}+2\Delta_2(a).\label{eq:gtmDsx2}
\end{align}
\end{lemma}

\section{The $D$- and $D_s$-optimal designs}

Before characterizing the $D$- and $D_s$-optimal designs, we first characterize a small class of designs guaranteed to contain $D$- and $D_s$-optimal designs, and then we obtain the optimal designs from the class. Such an approach is based on the so-called ``essentially complete classes'' used in \cite{Cheng95} and \cite{Yang09}, among others, to identify optimal designs in certain settings. Recently, a series of references, including \cite{Yang12} and \cite{Yang15}, proposed a framework to obtain essentially complete classes, and provided several important properties of optimal designs identified in the class. We follow such an approach.
A class $\mathcal{C}$ is called an essentially complete class if, for an arbitrary design $\xi_1\notin\mathcal{C}$, there exists a design $\xi_2\in\mathcal{C}$ such that $M(\xi_2)-M(\xi_1)$ is nonnegative definite. Therefore, it guarantees that there exists an optimal design in $\mathcal{C}$ when the design criterion has certain properties such as concavity, isotonicity, and smoothness.  For example, Assumption A in \cite{Yang15} includes criteria like $D$- and $D_s$-optimality, and other criteria based on the information matrix. Let 
\begin{align*}
\mathcal{C}_0=\left\{\{(x_i,w_i)\}_{i=1}^3:~x_L=x_1<x_2<x_3\leq x_U, \sum_{i=1}^3w_i=1, w_i\geq0\right\}.
\end{align*}
By allowing $w_i=0$  at some points, this class contains all one- and two-point designs, and all three-point designs with a support point at  $x_L$. By applying Theorem 2(a) in \cite{Yang12}, we have the following theorem which allows us to search for a $D$- and a $D_s$-optimal group testing design in $\mathcal{C}_0$.
\begin{theorem}\label{thm:cclass}
The class $\mathcal{C}_0$ is essentially complete. 
\end{theorem}
The justifications of this result and others are given in the Appendix.

\subsection{The unique $D$-optimal design}

When we wish to estimate prevalence, sensitivity, and specificity simultaneously, we need a design with at least three points to estimate the three parameters. The sought design is a $D$-optimal design for estimating the three parameters of the model. From Theorem \ref{thm:cclass}, such a $D$-optimal group testing design exists in $\mathcal{C}_0$ and the following theorem ensures the uniqueness of the $D$-optimal design and characterizes it.
\begin{theorem}\label{thm:gtmD}
The $D$-optimal design $\xi_D$ for estimating the prevalence and the two testing error rates
is unique. It has three group sizes $x_L$, $x_2^*$, and $x_U$ with equal weights,
where \mbox{$x_2^*-x_L=\log(A_1)/\log(1-p_0)$} and $A_1$ uniquely solves (\ref{eq:gtmDx2}).
\end{theorem}
Theorem \ref{thm:gtmD} asserts that the $D$-optimal design requires equal numbers of trials at three distinct group sizes, with two of them being the extreme sizes $x_L$ and $x_U$. No design with four or more group sizes can estimate the three parameters more accurately than the $D$-optimal design. In addition, our proof of Theorem \ref{thm:gtmD} implies that (a) the root of equation (\ref{eq:gtmDx2}) must be unique in $(r,1)$,  or otherwise it contradicts the uniqueness of the $D$-optimal design, and (b) either having a smaller value of $x_L$ or having a larger value of $x_U$ always strictly reduces the generalized variance of the MLE of $\theta$. Thus, we suggest the lower bound $x_L$ to be one whenever possible.  However, too large a value of $x_U$ may cause a dilution effect \citep{Zenios98}, impacting the sensitivity or the specificity of the test.  Thus, the upper bound $x_U$ should be set carefully.
When $p_0$ and $x_L$ are small and $x_U$ is large, we observe that $\pi(x_L)$ is close to $1-p_2$ and $\pi(x_U)$ is close to $p_1$, which implies that trials at size $x_L$ and $x_U$ mainly provide information for $p_2$ and $p_1$ respectively. We may conclude that $x_2^*$ contributes most for estimating $p_0$. In practice, $p_1$ and $p_2$ for a given test may be estimated by comparing the test results with a gold standard test having no testing errors. In the absence of a gold standard, the maximum likelihood approach can be seen as comparing the test results from the intermediate group size with those from the two extreme group sizes, to acquire information about the error parameters $p_1$, $p_2$, and the prevalence $p_0$. Accordingly, we may regard the test results from the two extreme group sizes as playing a similar role to that of a gold standard test for providing information on $p_1$ and $p_2$.

\subsection{The unique $D_s$-optimal design}

Theorem \ref{thm:gtmD} assumes that there is equal interest in estimating prevalence, sensitivity and specificity. In practice, the interest is usually to accurately estimate the prevalence alone, and the sensitivity and specificity of the test are treated as unknown nuisance parameters. When the main goal is only to estimate the prevalence as accurately as possible, the $D_s$-optimality criterion (\ref{eq:Dscriterion}) is the most suitable. A $D_s$-optimal design minimizes the $(1,1)$ element of the inverse of the information matrix (\ref{eq:inforM2}), which is proportional to the asymptotic variance of the prevalence estimator.
To determine the optimal weights under the $D_s$-criterion for a fixed number of group sizes, suppose that there are three group sizes $x_1<x_2<x_3 $  in the allowable range $[x_L,x_U]$. Let
\begin{align*}
Q_i(x_1,x_2,x_3)=\lambda(x_i)^{-1}\left\{(1-p_0)^{x^{(i)}_1}-(1-p_0)^{x^{(i)}_2}\right\}^2>0,
\end{align*}
for $i=1,2,3$, where $x^{(i)}_1$ and $x^{(i)}_2$ are the two group sizes of $\{x_1,x_2,x_3\}\setminus\{x_i\}$ such that $x^{(i)}_1<x^{(i)}_2$. To be more specific, $x^{(1)}_1=x_2$, $x^{(1)}_2=x_3$, $x^{(2)}_1=x_1$, $x^{(2)}_2=x_3$, $x^{(3)}_1=x_1$, and $x^{(3)}_2=x_2$. For any group testing design with three group sizes, the next lemma prescribes the optimal weights under the $D_s$-criterion. 
\begin{lemma}\label{thm:gtmDsW}
For any group testing design with three group sizes, $x_1$, $x_2$, and $x_3\in[x_L,x_U]$, the optimal weights under  the $D_s$-criterion for only estimating the prevalence are unique and given by
\begin{align}
&w_i^s=\frac{Q_i(x_1,x_2,x_3)^{1/2}}{\sum_{j=1}^3Q_j(x_1,x_2,x_3)^{1/2}}
\label{eq:DsW}
\end{align}
for $i=1,2,$ and 3, and
\begin{align}
\max\Phi_s[M(\xi)]&=\log\left|M_f(\xi)\right|^2-2\log\left(\sum_{i=1}^3Q_i(x_1,x_2,x_3)^{1/2}\right).
\label{eq:DsWmin}
\end{align}
\end{lemma}
Our next result gives a complete description of the $D_s$-optimal design for estimating the prevalence when the testing error rates are unknown.
\begin{theorem}\label{thm:gtmDs}
The $D_s$-optimal design $\xi_s$ for estimating the prevalence alone is unique and requires three group sizes $x_L$, $x_2^s$, and $x_U$, where $x^s_2-x_L=\log(A_s)/\log(1-p_0)$ and $A_s$ uniquely solves (\ref{eq:gtmDsx2}). The proportions of groups with such sizes are given in (\ref{eq:DsW}) in Lemma \ref{thm:gtmDsW}.
\end{theorem}
We observe that (\ref{eq:gtmDsx2}) approximates (\ref{eq:gtmDx2}) when $\Delta_1$ and $\Delta_2$ are both close to zero, which may occur in practice when $\pi(x_L)$ tends to zero and $\pi(x_U)$ tends to one. Consequently, $A_1\approx A_s$, and the intermediate support point $x_2^s$ of the $D_s$-optimal design for estimating the prevalence is close to the intermediate support $x_2^*$ of the corresponding $D$-optimal design.
Theorem \ref{thm:gtmDs} shows that the $D_s$-optimal design for estimating prevalence with unknown testing error rates has similar properties to those of the $D$-optimal design for estimating the three parameters. In particular, (a) no design with four or more different group sizes can estimate the prevalence more accurately than the $D_s$-optimal design with three different group sizes, (b) the root of equation (\ref{eq:gtmDsx2}) must be unique in $(r,1)$, (c) either having a smaller value of $x_L$ or having a larger value of $x_U$ strictly improves the accuracy of prevalence estimation, and (d) groups with size $x_L$ and $x_U$ mainly contribute information about the specificity and sensitivity, respectively, and the groups with the intermediate size mainly contribute information about the prevalence.
When $x_L$ is small and $x_U$ is large,  which implies that $\pi(x_L)$ is much smaller than $0.5$ and $\pi(x_U)$ is much larger than $0.5$ by (\ref{eq:gt_resp}), then $Q_2(x_L,x_2^s,x_U)$ is much larger than $Q_1(x_L,x_2^s,x_U)$ and $Q_3(x_L,x_2^s,x_U)$. Consequently, the proportion of groups with the intermediate size is much larger than the other two group sizes. This observation is not surprising partly because the $D_s$-optimal design criterion seeks information for the prevalence only and does so by having more groups with the intermediate size.

\section{Design performance}

The optimal designs discussed in the previous section are locally optimal approximate designs, and their group sizes and the number of trials at each group size must be rounded to positive integers for implementation. In addition,  their information matrices are based on the asymptotic distribution of the MLE of the parameter vector. In what follows, we study the properties of the rounded optimal designs when we have moderate sample size (total number of trials), and their sensitive to the prespecified value of the parameter vector.
Recall that $\theta=(p_0,p_1,p_2)^{\mathrm{T}}$ is the true but unknown vector of parameters in our model, $[x_L,x_U]$ is the user-specified allowable range of group sizes, and we are interested in either estimating $\theta$ or only the prevalence parameter $p_0$. For an approximate design $\xi=\{(x_i,w_i)\}_{i=1}^k$ with integer-valued group sizes, denote its exact design with $n$ trials as $\xi{(n)}=\{(x_i,n_i/n)\}_{i=1}^k$, where $n_1,\ldots,n_k$ are integers obtained by applying the efficient rounding apportionment \citep[section 12.5]{puk06} to $\{nw_1,\ldots,nw_k\}$. Let $\hat{\theta}_{\xi(n)}$ be the MLE of the parameter vector $\theta$ under the exact design $\xi(n)$. We define the mean squared error matrix of  $\hat{\theta}_{\xi(n)}$ (scaled by sample size $n$) by
\begin{align}
\MSE(\hat{\theta}_{\xi(n)})=n\mathrm{E}\left[(\hat{\theta}_{\xi(n)}-\theta)(\hat{\theta}_{\xi(n)}-\theta)^{\mathrm{T}}\right].
\end{align}
Note that $\MSE(\hat{\theta}_{\xi(n)})$ converges to $M(\xi)^{-1}$ as $n$ goes to infinity by the asymptotic property of MLE.
The analytical form of $\MSE(\hat{\theta}_{\xi(n)})$ is quite complicated unless the sample size $n$ is very small. We therefore study $\MSE(\hat{\theta}_{\xi(n)})$ by simulation as follows. Let $N$ be the number of simulation replications. For $t=1,\ldots,N$, we simulate a sample of $n$ trials generated from $\xi(n)$, consisting of $n_i$ binary outcomes with response probability $\pi(x_i|\theta)$ for $i=1\ldots k$, and compute $\hat{\theta}_{\xi(n)}^{(t)}$, the MLE of $\theta$ from the sample. The simulation-based mean squared error matrix of $\xi$ is defined by 
\begin{align}
\widehat{\MSE}(\hat{\theta}_{\xi(n)})=\frac{n}{N}\sum_{t=1}^N(\hat{\theta}_{\xi(n)}^{(t)}-\theta)(\hat{\theta}_{\xi(n)}^{(t)}-\theta)^{\mathrm{T}}.
\label{eq:MSEn}
\end{align}
We use two measures to assess the finite sample performance of $\xi(n)$. They are the (simulation-based) $D$- and $D_s$-efficiencies, defined  respectively, by
\begin{align}
\EFF_{D}\{\xi(n)\}=\left\{\frac{|M(\xi_D)^{-1}|}{|\widehat{\MSE}(\hat{\theta}_{\xi(n)})|}\right\}^{1/3}
\quad\text{and}\quad
\EFF_{s}\{\xi(n)\}=\frac{M(\xi_s)^{-1}_{11}}{\widehat{\MSE}(\hat{\theta}_{\xi(n)})_{11}},
\label{eq:eff}
\end{align}
where $\xi_D$ and $\xi_s$ are the $D$- and $D_s$-optimal designs obtained from Theorems \ref{thm:gtmD} and \ref{thm:gtmDs}.
In the following, we set the number of simulation replications to $N=10000$, and set the prespecified values of the test parameters to match the Chlamydia study from Nebraska, the U.S.A. \citep[Table 1]{McMahan12}. \cite{McMahan12} prespecified the prevalence to be $0.07$, and the sensitivity and specificity to be  $0.93$ and $0.96$ for female's swab specimens. Suppose we have a budget for $n=3000$ trials, and the allowable range of the group size is $[x_L,x_U]=[1,61]$.

\subsection{Finite sample performance of the rounded optimal designs}

First, we construct the $D$- and $D_s$-optimal designs, and show that the effects of rounding (both the group sizes and their number of support points) and the use of asymptotic approximations are small. Under $\theta=(0.07,0.93,0.96)^{\rm T}$ and $[x_L,x_U]=[1,61]$, the $D$-optimal design $\xi_D$ from Theorem \ref{thm:gtmD} is equally supported on $\{1,16.79,61\}$. For $n=3000$, its rounded design $\xi_D'(3000)$ has group sizes $\{1,17,61\}$ with 1000 groups per size.
Similarly, the $D_s$-optimal design $\xi_s$ obtained from Theorem \ref{thm:gtmDs} is supported on $\{1,15.68,61\}$. Its rounded design $\xi_s'(3000)$ is therefore supported on $\{1,16,61\}$, and has $\{393,1884,723\}$ trials at $\{1,16,61\}$ respectively. These numbers of trials are based on the $D_s$-optimal weights on $\{1,16,61\}$ by Lemma \ref{thm:gtmDsW}, which are slightly different from the weights of $\xi_s$. The simulation results in Table \ref{tab:DDsUdesigns} show that the $D$-efficiency of $\xi_D'(3000)$ and the $D_s$-efficiency of $\xi_s'(3000)$ are both close to one. This indicates that the performance of the rounded $D$- (and $D_s$-) optimal design is similar to the asymptotic performance of the  $D$- (and $D_s$-) optimal approximate design.

\begin{table}
\caption{The $D$- and $D_s$-efficiencies of the rounded $D$- and $D_s$-optimal designs, 
and the $k$-point uniform designs for $k=3,\ldots,6$.
\label{tab:DDsUdesigns}}
\begin{tabular}{|c|cc|cccc|}
\hline
Design & $\xi_D'(3000)$ & $\xi_s'(3000)$   & $\xi_U^{(3)}(3000)$ & $\xi_U^{(4)}(3000)$ & $\xi_U^{(5)}(3000)$ & $\xi_U^{(6)}(3000)$\\
\hline
$\EFF_{D}$ & 0.9904 & 0.8193               & 0.7988           & 0.8456          & 0.8133           & 0.7805\\
$\EFF_{s}$  & 0.6977 & 1.0052               & 0.3170           & 0.4851          & 0.5274           & 0.5509\\
\hline
\end{tabular}
\end{table}

Next, we compare the finite sample performance of the rounded $D$- and $D_s$-optimal designs to $k$-point uniform designs $\xi_U^{(k)}(n)$, supported on the positive integers $\{1+\frac{60i}{k-1}:~i=0,\ldots,k-1\}$, for $k=3,\ldots,6$. For example, $\xi_U^{(4)}(3000)$ has 750 trials at each group size of $\{1,21,41,61\}$. Table \ref{tab:DDsUdesigns} also shows the performance of these $k$-point uniform designs. We observe that, although the optimal designs $\xi_D$ and $\xi_s$ have exactly three support points, $\xi_U^{(3)}(3000)$ is not the best among these uniform designs regarding $D$- or $D_s$-efficiency. In contrast, $\xi_U^{(4)}(3000)$ performs best for estimating $\theta$ among these uniform designs, but its $D$-efficiency is about $0.85$ and is less than $\EFF_D(\xi_D'(3000))$. Furthermore, the design $\xi_U^{(6)}(3000)$ has highest $D_s$-efficiency among these uniform designs, but its $D_s$-efficiency is only about $0.55$ and is significantly less than $\EFF_s(\xi_s'(3000))$. This indicates that, when focusing on estimating $\theta$ (or $p_0$), the rounded $D$- (or $D_s$-) optimal design under the true parameter vector performs much better than the uniform designs above.

\subsection{Robustness of the rounded optimal designs}

\begin{figure}[!h]
\centering%
\subfigure{
\includegraphics[width=0.2\textwidth]{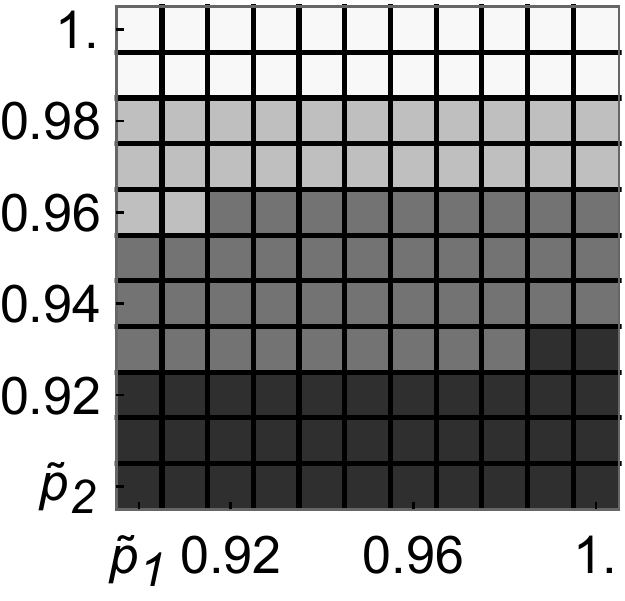}
}\hspace{.2cm}
\subfigure{
\includegraphics[width=0.2\textwidth]{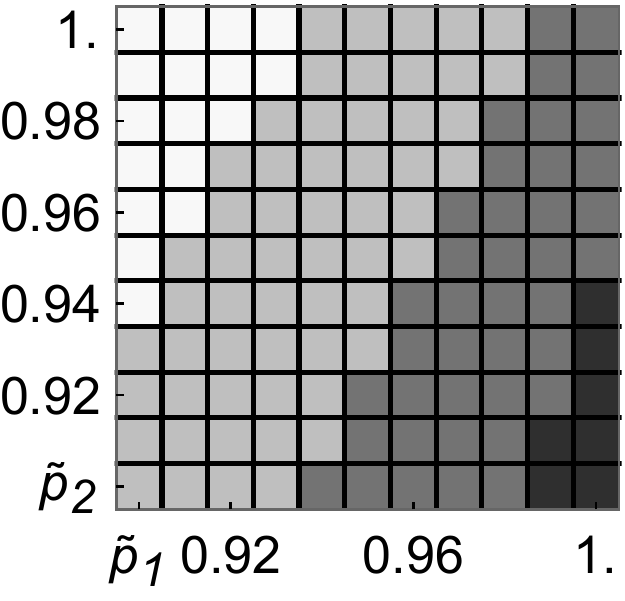}
}\hspace{.2cm}
\subfigure{
\includegraphics[width=0.2\textwidth]{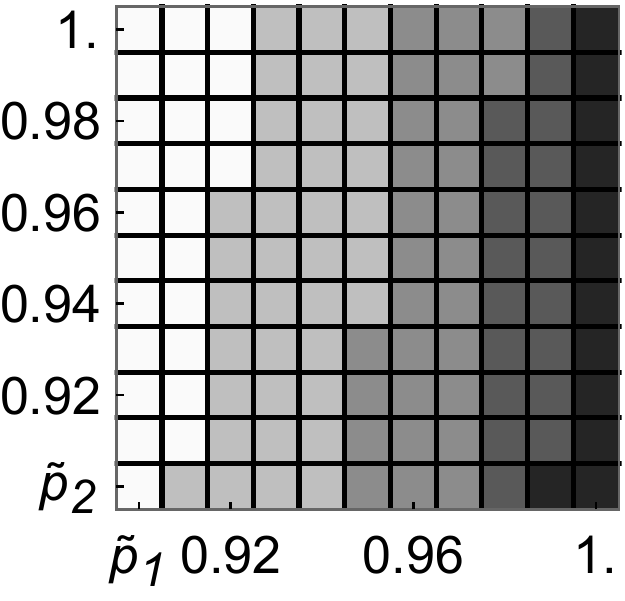}
}\hspace{.2cm}
\subfigure{
\includegraphics[width=0.2\textwidth]{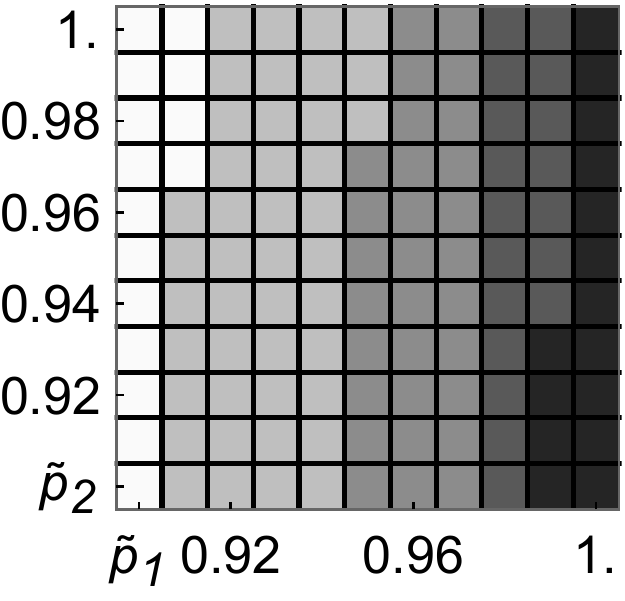}
}\\
\setcounter{subfigure}{0}
\vspace{-.2cm}
\hspace{.5cm}
\subfigure[$\tilde{p}_0=0.01$.]{
\includegraphics[width=0.15\textwidth]{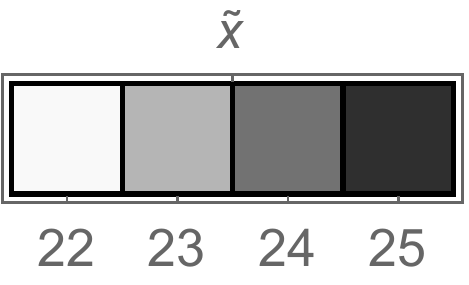}
}\hspace{1cm}
\subfigure[$\tilde{p}_0=0.04$.]{
\includegraphics[width=0.15\textwidth]{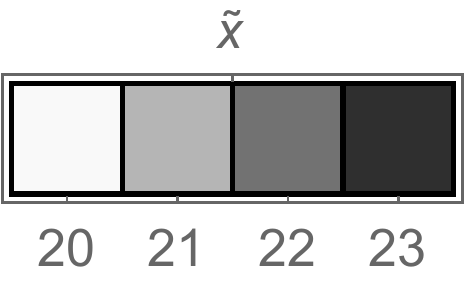}
}\hspace{1cm}
\subfigure[$\tilde{p}_0=0.07$.]{
\includegraphics[width=0.15\textwidth]{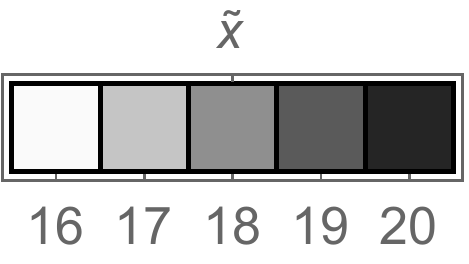}
}\hspace{1cm}
\subfigure[$\tilde{p}_0=0.10$.]{
\includegraphics[width=0.15\textwidth]{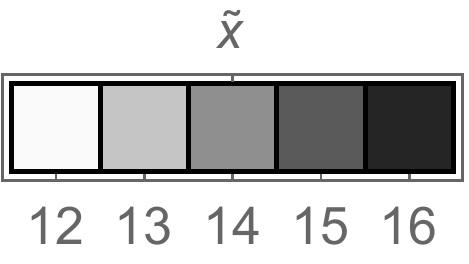}
}
\vspace{-.2cm}
\caption{The intermediate group size $\tilde{x}$ of the rounded $D$-optimal designs under a prespecified \mbox{$\tilde{\theta}=(\tilde{p}_0,\tilde{p}_1,\tilde{p}_2)^{\rm T}$}. Here $\tilde{p}_0=0.01, 0.04, 0.07$, and $0.10$; $\tilde{p}_1$ and $\tilde{p}_2=0.90, 0.91,\ldots, 0.99$, and $1.00$.}
\label{fig:x2D1}
\end{figure}

\begin{figure}
\centering%
\includegraphics[width=0.45\textwidth]{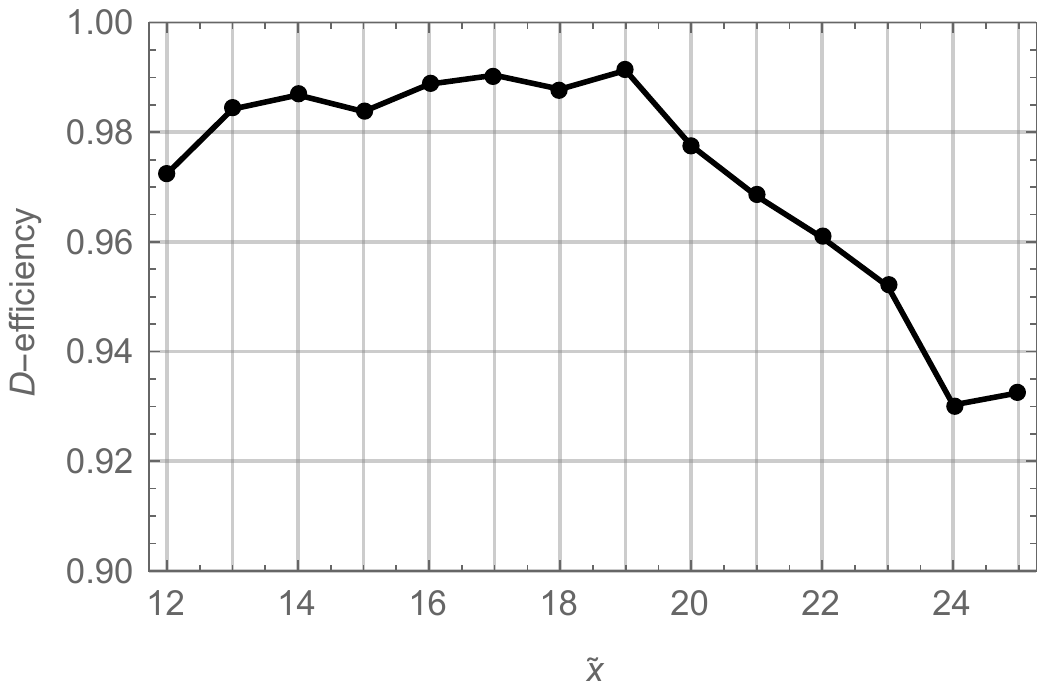}
\caption{The $D$-efficiency of the design that is equally supported on $\{1,\tilde{x},61\}$ with $n=3000$, for $\tilde{x}=12$ to $25$. 
These designs arise as rounded $D$-optimal designs for $\tilde{\theta}\in\Theta$ as defined in the text.
}\label{fig:DEff}
\end{figure}

\begin{figure}
\centering%
\subfigure[$\tilde{p}_0=0.01$.]{
\includegraphics[height=0.15\textheight]{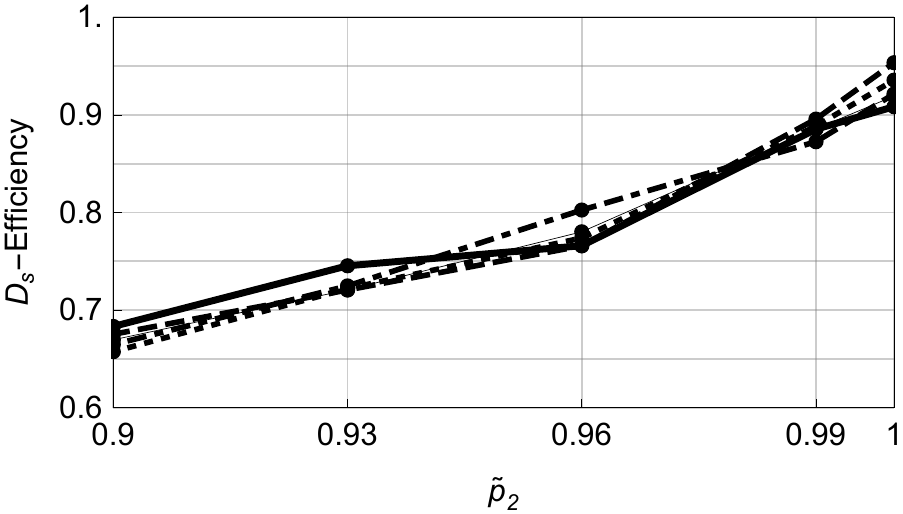}
}
\subfigure[$\tilde{p}_0=0.04$.]{
\includegraphics[height=0.15\textheight]{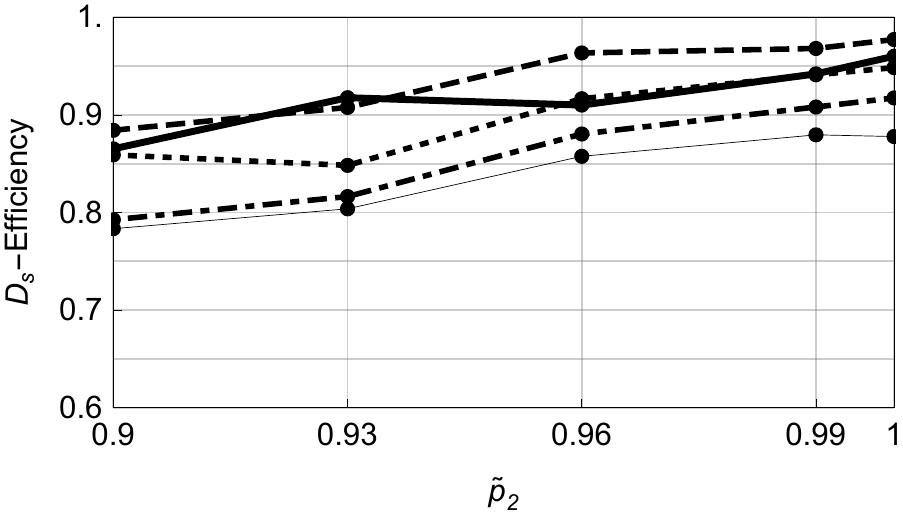}
}
\subfigure[$\tilde{p}_0=0.07$.]{
\includegraphics[height=0.15\textheight]{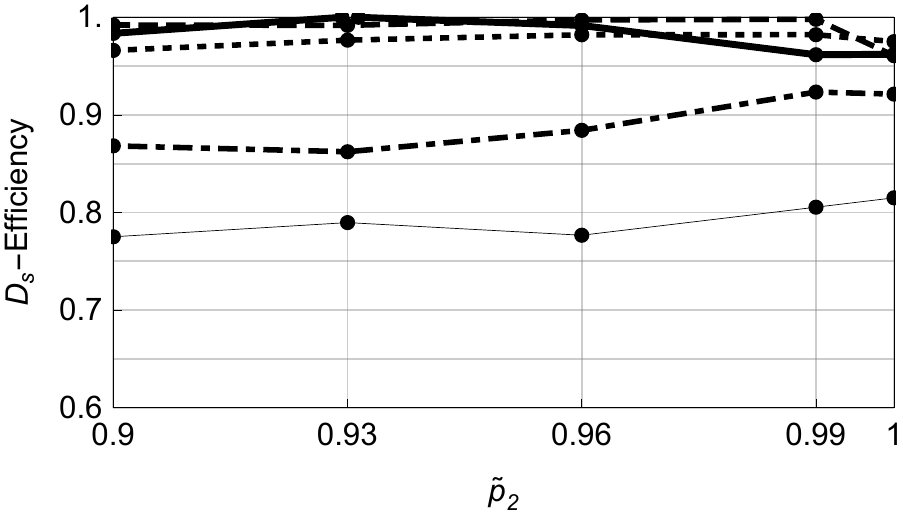}
}
\subfigure[$\tilde{p}_0=0.10$.]{
\includegraphics[height=0.15\textheight]{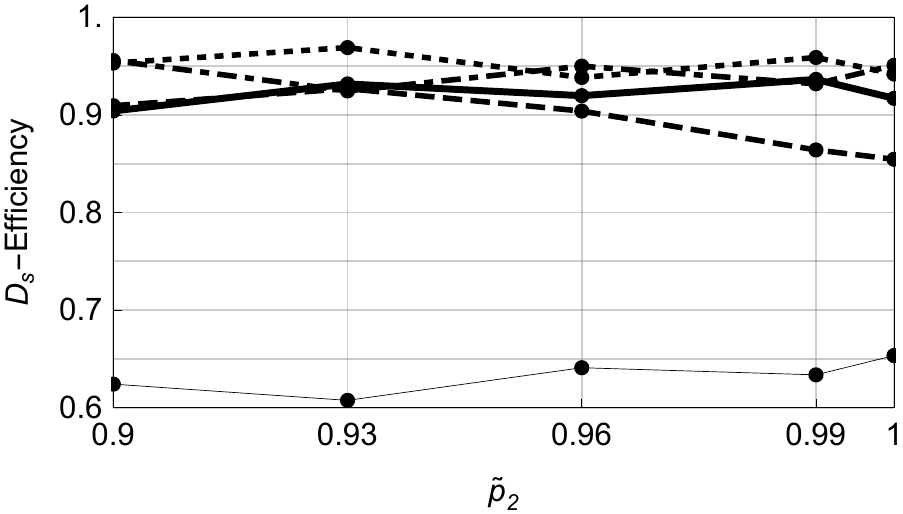}
}
\caption{The $D_s$-efficiency of the rounded $D_s$-optimal design $\xi_{s,\tilde{\theta}}'(3000)$ under prespecified $\tilde{\theta}=(\tilde{p}_0,\tilde{p}_1,\tilde{p}_2)^{\rm T}$. Here $\tilde{p}_0=0.01,0.04,0.07$, and $0.10$; $\tilde{p}_1=0.90$ (dashed), $0.93$ (thick), $0.96$ (dotted), $0.99$ (dot-dashed), and $1.00$ (thin); $\tilde{p}_2=0.90,~0.93,~0.96,~0.99$, and $1.00$. }
\label{fig:DsEff}
\end{figure}

Now we consider the situation where the prespecified value of $\theta$ used in the local design construction is incorrect. Specifically, suppose that the prespecified parameter vector $\tilde{\theta}=(\tilde{p}_0,\tilde{p}_1,\tilde{p}_2)^{\mathrm{T}}$ is a point in the region $\Theta=[0.01,0.10]\times[0.9,1]^2$, containing the true value of the parameter vector $\theta=(0.07,0.93,0.96)^{\mathrm{T}}$.
Let $\xi_{D,\tilde{\theta}}'(3000)$ be the rounded $D$-optimal design under the prespecified parameter vector $\tilde{\theta}$. In $D$-optimality, the prespecified value of the parameter vector only affects the intermediate group size $x_2^*$ in Theorem \ref{thm:gtmD}, but does not affect the other two group sizes and the weights on the three support points. Hence, $\xi_{D,\tilde{\theta}}'(3000)$ is equally supported on $\{1,\tilde{x},61\}$, where $\tilde{x}$ depends on $\tilde{\theta}$. Figure \ref{fig:x2D1} shows the corresponding $\tilde{x}$ for selected $\tilde{\theta}\in\Theta$. We observe that when the specified $\tilde{\theta}$ is within $\Theta$, $\tilde{x}$ shifts from $12$ to $25$, where $\tilde{x}=17$ for $\tilde{\theta}=\theta$; $\tilde{x}$ increases as $\tilde{p}_0$ or $\tilde{p}_2$ decreases, or  as $\tilde{p}_1$ increases. When the prevalence is small such as 0.01, the specificity seems to play a dominant role on the $D$-optimal design, whereas, as the prevalence gets larger, the dominant factor shifts from the specificity to sensitivity; and the sensitivity and specificity have similar effect on the design when the prevalence is close to $0.04$.
Figure \ref{fig:DEff} shows the $D$-efficiency of $\xi_{D,\tilde{\theta}}'(3000)$ for $\tilde{\theta}\in\Theta$ in terms of $\tilde{x}$, where $\tilde{x}=12$ to $25$. Generally, these rounded $D$-optimal designs under misspecified $\tilde{\theta}\in\Theta$ have $D$-efficiencies greater than $90\%$. Note that the rounded $D$-optimal design is equally supported on three group sizes including two extreme ones. In contrast, the design $\xi_U^{(3)}$ has the same form but its $D$-efficiency is only about $0.8$. This indicates that the rounded $D$-optimal design under a misspecified parameter vector is quite robust for estimating $\theta$.
Similar to the parameter misspecification setting for the $D$-optimal design, we obtain the rounded $D_s$-optimal designs designs $\xi_{s,\tilde{\theta}}'(3000)$ supported on $\{1,\tilde{x}_s,61\}$ under a possibly misspecified parameter vector $\tilde{\theta}$, where $\tilde{x}_s$ depends on $\tilde{\theta}\in\Theta$. A direct calculation shows that the intermediate group size $\tilde{x}_s$ is between 11 and 25, and the weights of group sizes 1, $\tilde{x}_s$, and 61 are respectively within $0.14\pm0.05$, $0.67\pm0.14$, and $0.19\pm0.12$.
Figure \ref{fig:DsEff} shows the $D_s$-efficiencies of $\xi_{s,\tilde{\theta}}'(3000)$ for selected $\tilde{\theta}\in\Theta$. We observe that $\EFF\{\xi_{s,\tilde{\theta}}'(3000)\}$  is larger than or close to $0.8$ when the prespecified $\tilde{\theta}$ belongs in the region $[0.04,0.10]\times[0.90,0.99]\times[0.90,1.00]\subset\Theta$, which limits the distance from $\tilde{\theta}$ to the true value of $\theta=(0.07,0.93,0.96)^{\rm T}$. However, when $\tilde{p}_0=0.01$ (not near the true prevalence $p_0=0.07$) or $\tilde{p}_1=1.00$ (not near the true sensitivity $p_1=0.93$), $\EFF\{\xi_{s,\tilde{\theta}}'(3000)\}$ drops to as low as $0.6$. We observe that the $D_s$-optimal design is slightly more sensitive than the $D$-optimal design when there is substantial parameter misspecification, but still performs well when the prespecified parameter vector is not far from its true value.
Finally, comparing Figures \ref{fig:DEff} and \ref{fig:DsEff} to Table \ref{tab:DDsUdesigns}, we observe that when the prespecified values of parameters are moderately misspecified, a rounded $D$- (or $D_s$-) optimal design still has at least 0.9 $D$-efficiency (or 0.8 $D_s$-efficiency), but those uniform designs only have at most $0.85$ $D$-efficiency (or 0.55 $D_s$-efficiency). Our conclusion is that when we wish to estimate $\theta$ (or $p_0$), we would recommend implementing a rounded $D$- (or $D_s$-) optimal design over any of the uniform design that we had investigated.

\section{Discussion}

Our work addresses several design issues for group testing studies using optimal design theory. Under uncertain testing error rates, we provide locally optimal designs for estimating prevalence alone or jointly with sensitivity and specificity. Our designs are justified by rigorous theoretical analysis and they provide the most accurate parameter estimates among all group testing designs. By contrasting the $D$- and $D_s$-optimal designs, we observe that relatively fewer trials at extreme group sizes are required in the $D_s$-optimal design, where the sensitivity and specificity are treated as nuisance parameters in the $D_s$-criterion. In the optimal approximate design, the theoretical group sizes and the frequencies of the group sizes are not guaranteed to be integers. We show by simulations that the rounded $D$- and $D_s$-optimal designs still have very high efficiencies. In addition, our simulation results show that our designs are still quite robust when the prespecified values of parameters are not far from their true values.
Although our theoretical results apply in general, our simulation studies focused on settings with 3000 trials, and had a maximum group size of approximately 60. This can be compared to two published studies using group testing: one in Russia estimating Chlamydia prevalence \citep{Shipitsyna07}, which had three designs, involving 150 to 1500 trials, and one in the U.S.A. focusing on Chlamydia diagnosis \citep{McMahan12}, which used 7000 trials. Thus the number of trials in our simulation study is similar to what has been used in practice, where the estimation burden is lighter because the sensitivity and specificity are treated as known there. We also note that the number of individual samples in our simulation example is larger than what were used in these two real-world studies. Since our designs are optimal, this reflects the unavoidable cost of estimating the prevalence as well as sensitivity and specificity from a single dataset with no prior knowledge about any of the parameters.
Our locally optimal designs can be used in the multi-stage adaptive approach as in \cite{HOS94}, where they assumed no testing errors. For example, a $D$- or $D_s$-optimal design may be implemented based on any available prior knowledge from the first stage. In the subsequent stages, we construct a $D$- or $D_s$-optimal design based on the information about the parameters obtained in the previous stages. Alternatively, a Bayesian optimal design approach or a minimax approach \citep{Dette14} can be implemented. For example, the experimenters may obtain a design that maximizes $D_s$-optimality criterion averaged over the parameters with respect to their prior distribution, or a design that minimizes the largest possible variance of the prevalence estimator among different values for the parameters. Such Bayesian or minimax optimization problems are difficult to address theoretically and even numerically, but are of interest for future studies.
The results presented here are appropriate for a setting in which the testing errors occur randomly and are unrelated to the group size within a prespecified range. In practice, the plausible range of group sizes comes from cost considerations or the physical constraints imposed on the study. A careful decision about the range of values for the group sizes is important because a group size that is too large may cause dilution effects, which in turn, may reduce the sensitivity and the specificity of the test. Our group testing framework can potentially be extended to accommodate dilution effects and other complicating issues, such as stratified populations. For example, we may wish to estimate the  prevalence of a rare trait in multiple geographic regions with potentially different disease prevalences, or when the sensitivity and specificity of the test vary among the sub-populations. Such issues may be interesting directions for future research.
\section*{Acknowledgements}
All the authors thank the reviewers for their helpful and thoughtful comments on earlier versions
of this paper. The first two authors were supported by the Ministry of Science and Technology,
Taiwan, under grant MOST101-2118-M-110-002-MY2. The third author was supported in part
by the National Center of Theoretical Science, Taiwan. The research of Wong that is reported in
this publication was partially supported by the National Institute of General Medical Sciences
of the National Institutes of Health under award R01GM107639. The content is solely the
responsibility of the authors and does not necessarily represent the official views of the National
Institutes of Health.

\appendix

\section*{Appendix}
Before giving the proofs of Theorem \ref{thm:cclass}, Theorem \ref{thm:gtmD}, Lemma \ref{thm:gtmDsW}, and Theorem \ref{thm:gtmDs},
we introduce additional notation. Let
\begin{align*}
a_i&=(1-p_0)^{x_i-x_L}\in[r,1], ~\textrm{for }i=1,2,\ldots,k,\\
C(a)&=
\frac{1}{(c-a)(\delta c+a)}
{\small\left[\begin{array}{c}a\log(a)\\c-a\\ \delta c+a\end{array}\right]}
{\small\left[\begin{array}{c}a\log(a)\\c-a\\ \delta c+a\end{array}\right]^{\mathrm{T}}}
=
{\small\left[\begin{array}{ccc}
\frac{(a\log(a))^2}{(c-a)(\delta c+a)} & \frac{a\log(a)}{c-a} & \frac{a\log(a)}{\delta c+a}\\
\frac{a\log(a)}{c-a} & \frac{c-a}{\delta c+a} & 1\\
\frac{a\log(a)}{\delta c+a} & 1 & \frac{\delta c+a}{c-a}
\end{array}\right]}\nonumber\\
&=
{\small\left[\begin{array}{ccc}
\Psi_5(a) & \Psi_4(a) & \Psi_3(a)\\
\Psi_4(a) & \Psi_1(a) & \Psi_0\\
\Psi_3(a) & \Psi_0 & \Psi_2(a)
\end{array}\right]},\\
P_1&={\small\left[\begin{array}{ccc}
\frac{(p_1+p_2-1)(1-p_0)^{x_L-1}}{\log(1-p_0)} & (p_1+p_2-1)x_L(1-p_0)^{x_L-1} & 0\\
0 & -(1-p_0)^{x_L} & 1\\
0 & -(1-p_0)^{x_L} & 0
\end{array}\right]},\quad\text{and}\\
P_2&={\small\left[\begin{array}{ccc}
1 & 0 & 0\\
0 & -1 & c\\
0 & 1 & \delta c
\end{array}\right]^{-1}}.
\end{align*}

\subsection*{Proof of Theorem \ref{thm:cclass}}

We use Theorem 2(a) in \cite{Yang12} to show that $\mathcal{C}_0$ is essentially complete. First we rewrite the information matrix (\ref{eq:inforM2}) as
\begin{align}
M(\xi)&=\frac{1}{(p_1+p_2-1)^2(1-p_0)^{2x_L}}
(P_1{\cdot}P_2)\left(\sum_{i=1}^kw_iC(a_i)\right)
(P_1{\cdot}P_2)^{\mathrm{T}}.
\label{eq:inforM2c}
\end{align}
Let $f_{i,1}(a)=d\Psi_i(a)/da$ for $i=1,\ldots,5$ and let
\begin{align*}
f_{i,j}(a)=\frac{d}{da}\left(\frac{f_{i,j-1}(a)}{f_{j-1,j-1}(a)}\right),\quad\text{for $j=2,\ldots,5$ and $i=j,\ldots,5$.}
\end{align*}
Then for $a\in[r,1]\subset(0,1]$, we have
\begin{align*}
f_{1,1}(a)&=-\frac{c+\delta c}{(\delta c+a)^2}<0; &
f_{2,2}(a)&=-\frac{2(c+\delta c)(\delta c+a)}{(c-a)^3}<0; &
f_{3,3}(a)&=-\frac{(c-a)^2(2a+c)}{2a^2(c+\delta c)^2}<0; \\
f_{4,4}(a)&=-\frac{2(c+\delta c)}{(2a+c)^2}<0; &
f_{5,5}(a)&=\frac{4a+2c}{a(c+\delta c)}>0.
\end{align*}
Therefore, $\prod_{i=1}^5f_{i,i}(a)=\frac{4}{a^3(c-a)(\delta c+a)}>0$ for $a\in[r,1]$. Hence, the group testing model satisfies the conditions for Theorem 2(a) in \cite{Yang12}, and therefore the designs having at most two points together with the designs having three points including $x_L$ (which coincides with $a=1$) form an essentially complete class. In addition, according to the proof of Theorem 2 in \cite{Yang12}, it can be shown that there exists $s_0,\ldots,s_5=\pm1$ such that 
\begin{align}
\{s_i\Psi_i\}_{i=0}^5\text{ is a Chebyshev system on $[r,1]$.}\label{eq:Tsys}
\end{align}
This property will be used in the proof of Theorem \ref{thm:gtmDs}.

\subsection*{Proof of Theorem \ref{thm:gtmD}}
\label{pf:gtmD}

There are 3 steps in the proof. In step (a) we show that the $D$-optimal design is unique and requires three group sizes with equal frequencies, including one being the smallest group size $x_L$. In step (b), we show that the upper bound $x_U$ is a support point of the $D$-optimal design. In step (c), we determine the intermediate group size.
\begin{itemize}
\item[(a)] By Theorem \ref{thm:cclass} and the fact that $\theta$ is only estimable under a design with at least 3 group sizes, there exists a $D$-optimal design $\xi_D$ belonging to $\mathcal{C}_0$ with exactly 3 group sizes and one of them is $x_L$. Because the number of support points of the $D$-optimal design $\xi_D$ is equal to the number of unknown parameters, $\xi_D$ must be equally supported (see for example, section 8.12 of \cite{puk06}). By Theorem 2.5 in \cite{Yang15} and the fact that $D$- and $\Phi_0$-criteria are equivalent, it follows that $\xi_D$ is the unique $D$-optimal design.
\item[(b)] Following (a), we denote the unique $D$-optimal design as $\xi_D=\{(x_i^*,\frac{1}{3})\}_{i=1}^3$, where $x_L=x_1^*<x_2^*<x_3^*\leq x_U$. Recall that $a_i=(1-p_0)^{x_i^*-x_L}$ ($x_1^*=x_L$ implies that $a_1=1$) and $r=(1-p_0)^{x_U-x_L}$. It suffices to show that $\frac{\partial\Phi_D[M(\xi_D)]}{\partial a_3}<0$ for all $a_3\in[r,a_2)$, i.e., the $D$-criterion is a strictly decreasing function of $a_3$. Thus the criterion is maximized by taking $a_3$ as small as possible, which implies $x_3^*=x_U$.
By (\ref{eq:inforM2c}), the information matrix of $\xi_D$ can be written as
\begin{align*}
M(\xi_D)&=\frac{1}{(p_1+p_2-1)^2(1-p_0)^{2x_L}}P_1\left(\sum_{i=1}^3\frac{1}{3}P_2C(a_i)P_2^{\mathrm{T}}\right)P_1^{\mathrm{T}},
\end{align*}
where
\begin{align*}
\sum_{i=1}^3\frac{1}{3}P_2C(a_i)P_2^{\mathrm{T}}
&=
\sum_{i=1}^3\frac{1}{3(c-a_i)(\delta c+a_i)}
{\small\left[\begin{array}{c}a_i\log(a_i)\\a_i\\1\end{array}\right]}{\small\left[\begin{array}{c}a_i\log(a_i)\\a_i\\1\end{array}\right]^{\mathrm{T}}}
=\Gamma_1\Gamma_2\Gamma_1^{\mathrm{T}};\\
\Gamma_1&={\small\left[\begin{array}{ccc}0 & a_2\log(a_2) & a_3\log(a_3)\\ 1 & a_2 & a_3 \\ 1 & 1 & 1\end{array}\right]};\\
\Gamma_2&={\small\left[\begin{array}{ccc}\frac{1}{3(c-1)(\delta c+1)} & 0 & 0 \\ 0 & \frac{1}{3(c-a_2)(\delta c+a_2)} & 0 \\ 0 & 0 & \frac{1}{3(c-a_3)(\delta c+a_3)}\end{array}\right]}.
\end{align*}
A direct calculation shows
\begin{align*}
\frac{\partial\Phi_D[M(\xi_D)]}{\partial a_3}
&=\frac{\partial}{\partial a_3}\log\left(\left|\Gamma_1\right|^2\left|\Gamma_2\right|\right)
=\frac{2}{\left|\Gamma_1\right|}\times\frac{\partial\left|\Gamma_1\right|}{\partial a_3}+\frac{1}{c-a_3}-\frac{1}{\delta c+a_3},
\end{align*}
and the following statements are true by calculus:
\begin{itemize}
\item[(i)] $\frac{1}{c-a_3}-\frac{1}{\delta c+a_3}<\frac{1}{c-a_3}-\frac{1}{c+a_3}=\frac{2a_3}{c^2-a_3^2}<\frac{2a_3}{1-a_3^2}$;
\item[(ii)] $\frac{\partial|\Gamma_1|}{\partial a_3}
    =1-a_2+a_2\log(a_2)+(1-a_2)\log(a_3)
    <1-a_2+a_2\log(a_2)+(1-a_2)\log(a_2)$\\
    $~\qquad=1-a_2+\log(a_2)<0$ for all $0<a_3<a_2<1$;
\item[(iii)] $|\Gamma_1|=a_3(1-a_2)\log(a_3)-a_2(1-a_3)\log(a_2)
    >\underset{a_3\rightarrow a_2}{\lim}|\Gamma_1|=0$;
\item[(iv)] $(1-a_3^2)\frac{\partial|\Gamma_1|}{\partial a_3}+a_3|\Gamma_1|
    =(1-a_3)(a_2\log(a_2)+a_3(1-a_2))-(1-a_2)(-\log(a_3)+(1-a_3)a_3)$<0.
\end{itemize}
The last assertion holds because $a_2\log(a_2)+a_3(1-a_2))<a_2\log(a_2)+a_2(1-a_2))<0$ and $(-\log(a_3)+(1-a_3)a_3)>0$ for all $0<a_3<a_2<1$. Consequently,
\begin{align*}
\frac{\partial\Phi_D[M(\xi_D)]}{\partial a_3}
<\frac{2}{\left|\Gamma_1\right|}
\times\frac{\partial\left|\Gamma_1\right|}{\partial a_3}+\frac{2a_3}{1-a_3^2}
=\frac{2}{\left|\Gamma_1\right|(1-a_3^2)}
\left((1-a_3^2)\frac{\partial\left|\Gamma_1\right|}{\partial{a_3}}
+a_3\left|\Gamma_1\right|\right)<0
\end{align*}
and the desired result holds.
\item[(c)] Our remaining goal is to determine $x_2^*$. We do this by showing that $\partial\Phi_D[M(\xi_D)]/\partial a_2=0$ is equivalent to (\ref{eq:gtmDx2}), and hence $x_2^*=x_L+\frac{\log A_1}{\log(1-p_0)}$, where $A_1$ is the solution of (\ref{eq:gtmDx2}). In (b) we have shown that the optimal choice for $a_3$ is $a_3^*=(1-p_0)^{x_U-x_L}=r$. By arguments similar to step (b), we have that
\begin{align*}
0
&=
\frac{\partial\Phi_D[M(\xi_D)]}{\partial{a_2}}
=
\frac{2}{\left|\Gamma_1\right|}\times
\frac{\partial\left|\Gamma_1\right|}{\partial a_2}
+\frac{1}{c-a_2}-\frac{1}{\delta c+a_2}\\
&=\frac{2\{(1-r)(1+\log(a_2))+r\log(r)\}}{(1-r)a_2\log(a_2)+a_2r\log(r)-r\log(r)}
+\frac{1}{c-a_2}-\frac{1}{\delta c+a_2}\\
&=\frac{2}{a_2}\left(1+\frac{1+\Delta_0\times\frac{1}{a_2}}
{\log(a_2)-\Delta_0\left(\frac{1}{a}-1\right)}\right)
+\frac{1}{c-a_2}-\frac{1}{\delta c+a_2}
\end{align*}
which is equivalent to (\ref{eq:gtmDx2}), and hence the proof is complete.
\end{itemize}

\subsection*{Proof of Lemma \ref{thm:gtmDsW}}

Let $\xi$ be a design supported on $\{x_1,x_2,x_3\}$ with $x_1<x_2<x_3$.  Our problem here is to find the vector of the positive weights  $\{w_i^s\}_{i=1}^3$ at these three given points that minimizes $\left(M(\xi)^{-}\right)_{11}$, where $M(\xi)=M_{f}(\xi)\cdot \textrm{Diag}_\lambda(\xi)\cdot M_{f}(\xi)^\mathrm{T}$, $M_{f}(\xi)=(f(x_1) ~ f(x_2) ~ f(x_3))$ is nonsingular and $\textrm{Diag}_\lambda(\xi)=diag(w_1^s\lambda(x_1),\ldots,w_3^s\lambda(x_3))$ is positive-definite. Let $e_1=(1,0,0)^\mathrm{T}$. Then
\begin{align*}
M_{f}(\xi)^{-1}\cdot e_1&=
\frac{1}{\left|M_{f}(\xi)\right|}
\left(\begin{array}{c}
(1-p_0)^{x_2}-(1-p_0)^{x_3}\\(1-p_0)^{x_3}-(1-p_0)^{x_1}\\(1-p_0)^{x_1}-(1-p_0)^{x_2}
\end{array}\right)
\end{align*}
and we have
\begin{align}
\left(M(\xi)^{-}\right)_{11}&=e_1^\mathrm{T}\cdot(M_{f}(\xi)^{-1})^{T}\cdot \textrm{Diag}_\lambda(\xi)^{-1}\cdot M_{f}(\xi)^{-1}\cdot e_1\nonumber\\
&=\left(M_{f}(\xi)^{-1}\cdot e_1\right)^\mathrm{T}\cdot \textrm{Diag}_\lambda(\xi)^{-1}\cdot\left(M_{f}(\xi)^{-1}\cdot e_1\right)\nonumber\\
&=\frac{1}{\left|M_{f}(\xi)\right|^2}\sum_{i=1}^3\frac{Q_i(x_1,x_2,x_3)}{w_i^s}.
\label{arg:ws}
\end{align}
Since $Q_i(x_1,x_2,x_3)>0$ for $i=1,2,3$ and $\left|M_{f}(\xi)\right|^2>0$, we apply the method of Lagrange multipliers directly to minimize the value in (\ref{arg:ws}) subject to the constraints on the weights.  The resulting solution is displayed in (\ref{eq:DsW}).
\subsection*{Proof of Theorem \ref{thm:gtmDs}}
We prove this theorem by similar steps as those in the proof of Theorem \ref{thm:gtmD}. In Step (a), we show that (i) $p_0$ is only estimable under a design with at least three distinct group sizes, and (ii) the $D_s$-optimal design is unique. Therefore, by Theorem \ref{thm:cclass}, the unique $D_s$-optimal design has exactly three group sizes, and one of them is the smallest allowable group size $x_L$. Steps (b) and (c) of this theorem's proof use similar arguments as those in the proof of Theorem \ref{thm:gtmD} and have therefore been omitted.
\begin{itemize}
\item[(i)]
The result is shown by contradiction. Without loss of generality, suppose there exists a design $\xi=\{(x_i,w_i)\}_{i=1}^2$ such that $p_0=e_1^{\mathrm{T}}\theta$ is estimable under $\xi$, where \mbox{$x_L\leq x_1<x_2\leq x_U$}, \mbox{$w_1,w_2\geq0$}, \mbox{$w_1+w_2=1$}, and $e_1=(1,0,0)^\mathrm{T}$. Therefore, $e_1$ belongs to the range of $M(\xi)=M_{f}(\xi)\cdot \textrm{Diag}_\lambda(\xi)\cdot M_{f}(\xi)^\mathrm{T}$, where $M_{f}(\xi)=[f(x_1) ~ f(x_2)]$ and $\textrm{Diag}_\lambda(\xi)=\text{diag}(w_1\lambda(x_1),w_2\lambda(x_2))$. Hence, $e_1$ belongs to the range of $M_{f}(\xi)$, or equivalently, the determinant of $(f(x_1) ~ f(x_2) ~ e_1)=0$. However,
\begin{align*}
\left|f(x_1) ~ f(x_2) ~ e_1\right|
&=
{\small\left|\begin{array}{ccc}
x_1(p_1+p_2-1)(1-p_0)^{x_1-1} & x_2(p_1+p_2-1)(1-p_0)^{x_2-1} & 1\\
1-(1-p_0)^{x_1} & 1-(1-p_0)^{x_2} & 0\\
-(1-p_0)^{x_1} & -(1-p_0)^{x_2} & 0
\end{array}\right|}\\
&=(1-p_0)^{x_1}-(1-p_0)^{x_2}>0
\end{align*}
for arbitrary $x_1<x_2$ and $p_0\in(0,1)$. This contradiction shows that $p_0$ is only estimable under a design with at least three points.
\item[(ii)] Suppose that there are two different $D_s$-optimal designs, $\xi_1$ and $\xi_2$, with at least 3 points. Let $\xi=\frac{1}{2}\xi_1+\frac{1}{2}\xi_2$. Therefore, $\xi$ must have at least 4 different support points due to $\xi_1\neq\xi_2$ and Lemma \ref{thm:gtmDsW}. By the concavity of $\Phi_s$, we have that $\Phi_s[M(\xi)]\geq\Phi_s[M(\xi_1)]=\Phi_s[M(\xi_2)]$, and hence $\xi$ is also $D_s$-optimal. By the equivalence theorem (see for example, section 2.7 of \cite{Fed72} or section 10.3 of \cite{atk07}), we must have
\begin{align*}
\phi_s(x,\xi)=1-\lambda(x)f(x)^\mathrm{T}M(\xi)^{-}f(x)+\lambda(x)f_s(x)^\mathrm{T}M_s(\xi)^{-}f_{s}(x)\geq0
\quad\text{ for all }x\in[x_L,x_U],
\end{align*}
with equality at each support point of $\xi$, where $f_s(x)=(f_1(x),f_2(x))^\mathrm{T}$ and $M_s(\xi)$ is the $2\times2$ submatrix of $M(\xi)$ deleting its first row and first column. Therefore, there exists a small enough $\epsilon>0$ such that $\phi_s(x,\xi)-\epsilon$ has at least $4\times2-2=6$ roots. In contrast, by direct calculation, $\phi_s(x,\xi)-\epsilon=\sum_{i=0}^5d_is_i\Psi_i(a)$ for some $d_i\in\mathbb{R}$ with $\sum_{i=0}^5d_i^2>0$, where $a=(1-p_0)^{x-x_L}\in[r,1]$ and $\{s_i\Psi_i\}_{i=0}^5$ is a Chebyshev system on $[r,1]$ by (\ref{eq:Tsys}). Therefore, Theorem 4.1 on page 22 of \cite{Karl66} shows that $\phi_s(x,\xi)-\epsilon$ has at most 5 roots. This contradiction shows that the $D_s$-optimal design is unique.
\end{itemize}

\end{document}